\documentclass[11pt]{article}
\usepackage{amsmath, amssymb}
\usepackage{latexsym}
\begin{document}
\newcommand{\qed}{\hfill \ensuremath{\square}}
\newtheorem{thm}{Theorem}[section]
\newtheorem{cor}[thm]{Corollary}
\newtheorem{lem}[thm]{Lemma}
\newtheorem{prop}[thm]{Proposition}
\newtheorem{defn}[thm]{Definition}
\newtheorem{rem}[thm]{Remark}
\newtheorem{ex}[thm]{Example}
\bibliographystyle{plain}
\numberwithin{equation}{section}

\def\nm{\noalign{\medskip}}
\newcommand{\Om}{\Omega}
\newcommand{\Real}{\mathbb{R}}
\newcommand{\nuu}{\tilde{\nu}}
\newcommand{\bohm}{{\partial}{\ohm}}
\newcommand{\la}{\langle}
\newcommand{\ra}{\rangle}
\newcommand{\ms}{\mathcal{S}_\ohm}
\newcommand{\mk}{\mathcal{K}_\ohm}
\newcommand{\mks}{\mathcal{K}_\ohm ^{\ast}}
\newcommand{\grad}{\bigtriangledown}
\newcommand{\ds}{\displaystyle}
\newcommand{\pf}{\medskip \noindent {\sl Proof}. ~ }
\newcommand{\p}{\partial}
\newcommand{\z}{\zeta}
\newcommand{\eps}{\varepsilon}
\newcommand{\pd}[2]{\frac {\p #1}{\p #2}}
\newcommand{\dbar}{\overline \p}
\newcommand{\eqnref}[1]{(\ref {#1})}
\newcommand{\na}{\nabla}
\newcommand{\ep}{\epsilon}
\newcommand{\vp}{\varphi}
\newcommand{\Scal}{\mathcal{S}}
\newcommand{\Ccal}{\mathcal{C}}
\newcommand{\Dcal}{\mathcal{D}}
\newcommand{\Lcal}{\mathcal{L}}
\newcommand{\Kcal}{\mathcal{K}}
\newcommand{\Ecal}{\mathcal{E}}
\newcommand{\Hcal}{\mathcal{H}}
\newcommand{\Ncal}{\mathcal{N}}
\newcommand{\Abar}{\overline A}
\newcommand{\Cbar}{\overline C}
\newcommand{\Ebar}{\overline E}
\newcommand{\RR}{\mathbb{R}}
\newcommand{\CC}{\mathbb{C}}
\newcommand{\GG}{\mathbb{G}}
\newcommand{\MM}{\mathbb{M}}
\newcommand{\II}{\mathbb{I}}
\newcommand{\tr}{\textrm{tr}\,}
\def\Bu {{\mathbf u}}
\def\Bv {{\mathbf v}}
\def\Bw {{\mathbf w}}
\def\Bx {{\mathbf x}}
\def\By {{\mathbf y}}
\def\Ba {{\mathbf a}}
\def\Bb {{\mathbf b}}
\def\Bf {{\mathbf f}}
\def\Bn {{\mathbf n}}
\def\BI {{\mathbf I}}
\def\BA {{\mathbf A}}
\def\BB {{\mathbf B}}
\def\BP {{\mathbf P}}
\def\BQ {{\mathbf Q}}
\def\BR {{\mathbf R}}
\def\BU {{\mathbf U}}
\def\BGam {{\mathbf{\Gamma}}}
\newcommand{\beq}{\begin{equation}}
\newcommand{\eeq}{\end{equation}}


\title{Progress on the Strong Eshelby's Conjecture and Extremal Structures for the Elastic Moment
Tensor\thanks{This work was partially supported by the ANR project
EchoScan (AN-06-Blan-0089), the STAR project 190117RD,
KRF-2008-220-C00002, and research grants at the Inha University, NSF
through the grant DMS-070978, and BK21 at KAIST. }}

\author{Habib Ammari\thanks{\footnotesize Laboratoire Ondes et Acoustique, CNRS UMR 7587, ESPCI,
10 rue Vauquelin, 75231 Paris Cedex 05, France
(habib.ammari@polytechnique.fr).} \and Yves
Capdeboscq\thanks{\footnotesize Mathematical Institute 24-29 St
Giles' Oxford OX1 3LB, UK (capdeboscq@maths.ox.ac.uk).} \and
Hyeonbae Kang\thanks{\footnotesize Department of Mathematics, Inha
University, Incheon 402-751, Korea (hbkang@inha.ac.kr, hdlee@inha.ac.kr).} \and Hyundae Lee\footnotemark[4]
\and Graeme W. Milton\thanks{\footnotesize Department of Mathematics, University
of Utah, Salt Lake City, UT 84112, USA (milton@math.utah.edu).} \and Habib Zribi\thanks{\footnotesize
Department of Mathematical Sciences, Korean Advanced Institute of
Science and Technology, Daejeon 305-701, Korea
(zribi@cmapx.polytechnique.fr). }}

\date{} \maketitle

\begin{abstract}
We make progress towards proving the strong Eshelby's conjecture in three
dimensions. We prove that if for a single nonzero uniform loading the
strain inside inclusion is constant and further the eigenvalues of this
strain are either all the same or all distinct, then the inclusion must be
of ellipsoidal shape. As a consequence, we show that for two linearly independent loadings
the strains inside the inclusions are uniform, then the inclusion must be
of ellipsoidal shape. We then use this result to address
a problem of determining the shape of an inclusion when the elastic moment
tensor (elastic polarizability tensor) is extremal. We show that the shape
of inclusions, for which the lower Hashin-Shtrikman bound either
on the bulk part or on the shear part of the elastic moment tensor is
attained, is an ellipse in two dimensions and an ellipsoid in three
dimensions.
\end{abstract}

%
\section{Introduction and statements of results}
%

In theory of composites or micro-structures, it is important to find
inclusion shapes which produce the minimal energy. In
relation to such shapes Eshelby \cite{esh57} showed that if the
inclusion is of ellipsoidal shape, then for any uniform loading the
strain inside $\Om$ is uniform. We call this remarkable property
Eshelby's uniformity property. Eshelby then conjectured in
\cite{esh61} that ellipsoids are the only shape (structure) with
such a uniformity property.

Eshelby's conjecture may be interpreted in two different ways:

\medskip

\noindent {\bf Weak Eshelby's conjecture}. If the strain is constant
inside $\Om$ for all loadings, then $\Om$ is an ellipse (2D) or an
ellipsoid (3D).

\medskip

\noindent {\bf Strong Eshelby's conjecture}. If the strain is
constant inside $\Om$ for a single loading, then $\Om$ is an ellipse
(2D) or an ellipsoid (3D).

\medskip

\noindent The strong Eshelby conjecture of course implies the weak
one.

The strong Eshelby's conjecture has been proved to be true in two
dimensions by Sendeckyj \cite{sen70} (see also \cite{KM_ARMA_08,
Liu_PRSA_08} for alternative proofs). However it is only recently
that the weak Eshelby conjecture was proved to be true in three
dimensions: by Kang-Milton \cite{KM_ARMA_08} and Liu
\cite{Liu_PRSA_08}. We refer to above mentioned papers (and
\cite{Kan_MM_09}) for comprehensive account of developments on the
Eshelby conjecture.

Regarding the strong Eshelby conjecture in three dimensions,
important pro--gress has been made by Liu. 
He showed in
\cite{Liu_PRSA_08} that the conductivity version of the strong
Eshelby conjecture fails to be true completely (see
\cite{Liu_PRSA_08} for a precise statements.) However, the strong
Eshelby's conjecture (for elasticity) has not been proved or
disproved. In this paper we consider the strong Eshelby conjecture.
Even though we are not able to resolve the conjecture completely, we
obtain results which is stronger than the weak version of Eshelby's
conjecture (and weaker than the strong version). We show that if the
strain inside inclusion is constant and in addition the eigenvalues
of the constant strain are either all the same or all distinct, then
the inclusion is of ellipsoidal shape. We then use this result to show that
for two linearly independent loadings
the strains inside the inclusions are uniform, then the inclusion must be
of ellipsoidal shape. It is worth emphasizing that the weak Eshelby's conjecture requires 6 linearly independent loadings while the strong Eshelby's conjecture does a single loadings.

In order to present results in more precise way let us introduce
some notation. Let $\Om$ be a bounded domain with a Lipschitz
boundary in $\RR^d$, $d=2,3$. The domain $\Om$ is occupied by a
homogeneous isotropic elastic material whose Lam\'e parameters are
$\tilde\lambda$ and $\tilde \mu$. We assume that the background
(the matrix) is also homogeneous and isotropic, and its Lam\'e
parameters are $\lambda$ and $ \mu$. Then the elasticity tensors for the
matrix and the inclusion can be written respectively as
 \beq
 \CC^{0}:=\lambda\mathbf{I}\otimes\mathbf{I}+2\mu \II
 \quad\mbox{and} \quad
 \CC^{1}:=\tilde{\lambda}\mathbf{I} \otimes \mathbf{I} + 2\tilde{\mu}\II,
 \eeq
where ${\mathbf{I}}$ is the $d\times d$ identity matrix (2-tensor)
and $\II$ is the identity 4-tensor. The the elasticity tensor for
$\RR^d$ in the presence of the inclusion $\Om$ is then given by
 \beq
 \CC_\Om:=(1-1_\Om) \CC^0+1_\Om \CC^1,
 \eeq
where $1_\Om$ is the indicator function of $\Om$.

Let $\kappa$ and $\tilde{\kappa}$ be bulk moduli of
$\mathbb{R}^{d}\setminus\overline{\Om}$ and $\Om$, respectively,
namely,
 $$
 \kappa = d\lambda+2\mu \quad\mbox{and}\quad \tilde{\kappa}=
 d\tilde{\lambda}+2\tilde{\mu}, \quad d=2,3.
 $$
It is always assumed that the strong convexity condition holds, {\it
i.e.},
 \begin{equation}\label{eq:trans-2}
 \mu>0,\quad \kappa >0,\quad\tilde{\mu}>0\quad\mbox{and}\quad \tilde{\kappa}>0\;.
 \end{equation}
We also assume that
 \[
 (\lambda-\tilde{\lambda})(\mu-\tilde{\mu})>0\;,
 \]
which implies that $\CC^{1}-\CC^{0}$ is either positive or negative
definite as an operator on the space $M_d^s$ of all $d \times d$
symmetric matrices.

We consider the following problem of the Lam\'e system of linear
elasticity: For a given non-zero symmetric $d \times d$ matrix $\BA$
\begin{equation}\label{main-eqn}
\left\{
\begin{array}{ll}
\nabla\cdot \CC_{\Om} \Ecal(\Bu)=0 \quad & \mbox{in } \RR^d \\
\nm \Bu(\Bx)- \BA \Bx=O(|\Bx|^{1-d}) \quad & \mbox{as
}|\Bx|\rightarrow \infty,
\end{array}
\right.
\end{equation}
where $\Ecal(\Bu)$ is the strain tensor, {\it i.e.},
 $$
 \Ecal(\Bu):= \frac{1}{2} (\nabla \Bu + \nabla \Bu^{T}) \quad (T \mbox{ for
 transpose}).
 $$
The matrix $\BA$ represents a uniform loading at infinity.

In this paper we prove the following improvements of the weak Eshelby
conjecture for the three dimensional elasticity.

\begin{thm} \label{mainthm}
Let $\Om$ be a simply connected bounded domain in $\RR^3$ with a
Lipschitz boundary. If the strain tensor $\Ecal(\Bu)$ of the
solution $\Bu$ to \eqnref{main-eqn} is constant in $\Om$ for a
nonzero symmetric matrix $\BA$ and $\Ecal(\Bu)$ within $\Om$ has eigenvalues
which are either all distinct or all the same, then $\Om$ is an
ellipsoid.
\end{thm}

\begin{thm} \label{mainthm2}
Let $\Om$ be a simply connected bounded domain in $\RR^3$ with a
Lipschitz boundary. If the strain tensors of
solutions to \eqnref{main-eqn} for two linearly independent $\BA$'s are constant in $\Om$, then $\Om$ is an
ellipsoid.
\end{thm}

The second main result of this paper is on the shape of the
inclusion whose elastic moment tensor (elastic polarizability
tensor) has an extremal property. In order to explain the second
result, we take the following definition of the Elastic Moment Tensor
(henceforth denoted as the EMT) \cite[Lemma
10.3]{AK_book_07}: Let $\BA$ be a $d \times d$ matrix and let
$\Bu_{\BA}$ be the solution to \eqnref{main-eqn} corresponding to
$\BA$. Then the EMT $\MM$ associated with the inclusion $\Om$ and
the elasticity tensors $\CC^0$ and $\CC^1$ is a 4-tensor defined by
 \beq \label{EMT_def}
 \MM \BA = \int_\Om (\CC^1 - \CC^0)
 \Ecal(\Bu_{\BA}) \, d \Bx.
 \eeq
The EMT may be defined in many different but equivalent ways. It is
worth noticing that if the strain $\Ecal(\Bu_{\BA})=\BB$ is constant
in $\Om$, then
 \beq \label{EMT2}
 \MM \BA = |\Om|(\CC^1 - \CC^0) \BB,
 \eeq
where $|\Om|$ denotes the volume of $\Om$.

The EMT enjoys several important properties. For example, it is
symmetric and positive-definite or negative-definite on the space
$M_d^s$ of $d\times d$ symmetric matrices, depending on the sign of
$\tilde\mu-\mu$. The notion of EMT is being used in variety of
contexts such as detection of small elastic inclusions for
non-destructive evaluation and medical imaging \cite{ACI_SIIS_08,
AGKL, AGKL08, AKNT_JE_02, AK_book_07, KKL_IP_03, KKL_IP_07} and
effective medium theory \cite{AK_book_07, AKL_IUM_06, milton}.

Let us introduce more notation in order to recall the optimal trace
bounds (the Hashin-Shtrikman bounds) for the EMT. Let
 \[
 \mathbf{\Lambda}_{1}:=\frac{1}{d}{\mathbf{I}}\otimes{\mathbf{I}},
 \quad{\mathbf{\Lambda}}_{2}:=\II-{\mathbf{\Lambda}}_{1}.
 \]
Then the elasticity tensor $\CC^{0}$ may be written as
 \[
 \CC^{0}= d \kappa
 \mathbf{\Lambda_{1}}+2\mu\mathbf{\Lambda_{2}},
 \]
and likewise for $\CC^1$. Since for any $d\times d$ symmetric matrix
$\BA$, ${\mathbf{I}}\otimes{\mathbf{I}}(\BA) =\tr (\BA){\mathbf{I}}$
and ${\II}(\BA)=\BA$, one can immediately see that
 \begin{equation}\label{eq:mbda11}
 {\mathbf{\Lambda}}_{1}{\mathbf{\Lambda}}_{1}={\mathbf{\Lambda}}_{1},\quad{\mathbf{\Lambda}}_{2}
 {\mathbf{\Lambda}}_{2}={\mathbf{\Lambda}}_{2},\quad{\mathbf{\Lambda}}_{1}{\mathbf{\Lambda}}_{2}=0.
 \end{equation}

We are now able to recall the optimal trace bounds for the EMT. For
$d=2,3$, let
 \begin{align}
 K_1 &:= \frac{1}{d (\tilde{\kappa}-\kappa)} \frac{d\tilde{\kappa}+2(d-1)\mu}{d\kappa+2(d-1)\mu} ,
 \label{eq:kone} \\
 K_2 &:= \frac{1}{2(\tilde{\mu}-\mu)} \left[ \frac{d^{2}+d-2}{2}+2\left(\tilde{\mu}-\mu\right)
 \left(\frac{d-1}{2\mu}+\frac{d-1}{d\kappa +2(d-1)\mu}\right)\right]. \label{eq:ktwo}
 \end{align}
The following trace bounds were obtained by Lipton
\cite{Lip_JMPS_93} (see also \cite{CK_AMO_08}): Suppose $|\Om|=1$ and let $\MM$ be the EMT associated with $\Om$,
then we have
 \begin{align}
 \tr \left(\mathbf{\Lambda}_{1}\MM^{-1}\mathbf{\Lambda}_{1}\right)
 & \leq K_1 ,  \label{eq:bd-itrace-1} \\
 \tr \left(\mathbf{\Lambda}_{2}\MM^{-1}\mathbf{\Lambda}_{2}\right)
 & \leq K_2, \label{eq:bd-itrace-2}
 \end{align}
provided that $\tilde{\kappa} -\kappa >0$. (If $\tilde{\kappa}
-\kappa <0$, the inequalities change the direction.) Since
$\mathbf{\Lambda}_{1}\MM^{-1}\mathbf{\Lambda}_{1}$ and
$\mathbf{\Lambda}_{2}\MM^{-1}\mathbf{\Lambda}_{2}$ are block
diagonal components for $\MM^{-1}$, one can see that
 $$
 \tr  \MM^{-1} = \tr (\mathbf{\Lambda}_{1}\MM^{-1}\mathbf{\Lambda}_{1}) + \tr (\mathbf{\Lambda}_{2}\MM^{-1}\mathbf{\Lambda}_{2}),
 $$
and hence
 \begin{equation}\label{eq:lowbound}
 \tr  \MM^{-1} \le K_1 + K_2.
 \end{equation}

Note that $\mathbf{\Lambda}_{1}\MM \mathbf{\Lambda}_{1}$ and
$\mathbf{\Lambda}_{2}\MM \mathbf{\Lambda}_{2}$ are the bulk and
shear parts of $\MM$, respectively. We also note that
\eqnref{eq:bd-itrace-1} and \eqnref{eq:bd-itrace-2} are lower bounds
for $\MM$ since they are upper bounds for $\MM^{-1}$. It is worth
emphasizing that upper bounds for $\MM$ are also derived in
\cite{Lip_JMPS_93}. In \cite{CK_AMO_08}, it is shown that inclusions
$\Om$ whose trace is close to the upper bound must be infinitely
thin. The upper and lower bounds for the EMT may also be derived as
a low volume fraction limit of the Hashin-Shtrikman bounds for the
effective moduli of the two phase composites, which was obtained by
Zhikov \cite{Zhi88, Zhi91} and Milton-Kohn \cite{MK88}. Benveniste
\cite{Ben87} obtained the upper and lower bounds of EMTs when
those EMTs happen to be isotropic. (See also \cite{MMM81}.)

In this paper we are interested in the shape of the inclusion whose
EMT satisfies the equality in either \eqnref{eq:bd-itrace-1} or
\eqnref{eq:bd-itrace-2}. This is an isoperimetric inequality for the
EMT. In this direction we prove the following theorem:
\begin{thm} \label{thm:main-thm}
Let $\Om$ be a simply connected bounded domain in $\RR^3$ with a
Lipschitz boundary. Suppose $|\Om|=1$ and let $\MM$ be the EMT associated with $\Om$. If
the equality holds in either \eqnref{eq:bd-itrace-1} or
\eqnref{eq:bd-itrace-2}, then $\Om$ is an ellipse in two dimensions
and an ellipsoid in three dimensions.
\end{thm}
We remark that optimal shapes for a cavity (hole) in two
dimension were investigated by Cherkaev {\it et al} \cite{CGMS98}
and Milton {\it et al} \cite{MSM03}.

The dimension of the space of symmetric 4-tensors in the three
dimensional space is 21, and hence the equalities
\eqnref{eq:bd-itrace-1} and \eqnref{eq:bd-itrace-2} are satisfied on
a 19 ($21-2$) dimensional surface in tensor space. However ellipsoid
geometries (with unit volume) only cover a 5 dimensional manifold
within that 19 dimensional space.

It is interesting to notice similarity of Theorem \ref{thm:main-thm}
to the P\'olya-Szeg\"o conjecture, which asserts that the inclusion
whose polarization tensor has the smallest trace is a disk or a
ball. The P\'olya-Szeg\"o conjecture was proved to be true by
Kang-Milton \cite{KM_CONM_06, KM_ARMA_08}. As for the
P\'olya-Szeg\"o conjecture, Theorem \ref{thm:main-thm}, which
concerns elasticity, will be proved using Eshelby's conjecture.

In order to prove Theorem \ref{thm:main-thm}, we will show that if
equality holds in \eqnref{eq:bd-itrace-1}, then the strain
tensor corresponding to a certain uniform loading $\BA$ (with a
special structure) is constant in $\Om$, while if the equality holds
in \eqnref{eq:bd-itrace-2}, then the strain tensors corresponding to
five (two in 2D) linearly independent uniform loadings are constant
in $\Om$ . Thus in two dimensions the strong Eshelby conjecture
immediately implies that the inclusion is an ellipse. However, in
three dimension, the weak Eshelby conjecture does not guarantee that
the inclusion is an ellipsoid. In order to apply the weak Eshelby
conjecture, we need to have equalities in both
\eqnref{eq:bd-itrace-1} and \eqnref{eq:bd-itrace-2}, or the equality
in the whole lower trace bound \eqnref{eq:lowbound}. But we are able
to show additionally that if the equality holds in
\eqnref{eq:bd-itrace-1} then the eigenvalues of the strain tensor
are all the same, and that if the equality holds in
\eqnref{eq:bd-itrace-2} then strains corresponding to five linearly
independent loadings are constant.
Thus thanks to Theorem \ref{mainthm} and Theorem \ref{mainthm2} we are able to conclude that
the inclusion is of ellipsoidal shape.

This paper is organized as follows. In section 2, we show that the
displacement vectors can be decomposed in a way similar to the
Helmholtz decomposition. This is done using the single layer
potential for the Lam\'e system. Theorem \ref{mainthm} is proved in
section 3, Theorem \ref{mainthm2} in section 4, and Theorem
\ref{thm:main-thm} in section 4. The appendix is for the proof of
Lemma \ref{discrim} which is used to prove Theorem \ref{mainthm2}.

%
\section{Single layer potential}
%

Let us first recall the notion of the single layer potential for the
Lam\'e operator $\Lcal_{\CC}\Bu:=\nabla \cdot \CC \Ecal(\Bu)$. The
Kelvin matrix ${\bf \Gamma} = ( \Gamma_{ij} )_{i,j=1}^3$ of the
fundamental solution  to the Lam{\'e} operator $\Lcal_{\CC}$ in
three dimensions is given by
 \beq \label{Kelvin}
 \Gamma_{ij} (\Bx) : = - \frac{\alpha_1}{4 \pi} \frac{\delta_{ij}}{|\Bx|} -  \frac{\alpha_2}{4 \pi}
 \frac{x_i x_j}{|\Bx|^3}, \quad  \Bx \neq 0\;,
 \eeq
where
 \beq \label{ab} \alpha_1= \frac{1}{2} \left ( \frac{1}{\mu} +
 \frac{1}{2\mu + \lambda} \right ) \quad\mbox{and}\quad \alpha_2=
 \frac{1}{2} \left ( \frac{1}{\mu} - \frac{1}{ 2 \mu + \lambda}
 \right )\;.
 \eeq
The single layer potential of the vector valued density function
$\Bf$ on $\p \Om$ associated with the Lam{\'e} parameters $(\lambda,
\mu)$ is defined by
 \beq
 \Scal_\Om [\Bf] (\Bx)  := \int_{\p \Om} {\bf \Gamma} (\Bx-\By)
 \Bf (\By)\, d \sigma (\By)\;, \quad \Bx \in \RR^3\; .
 \label{singlelayer}
 \eeq

Using the divergence theorem, we have
 \begin{align*}
 \Scal_\Om[\mathbf{f}] (\Bx) & = - \frac{\alpha_1}{4\pi} \int_{\p \Omega} \frac{\mathbf{f}(\By)}{|\Bx-\By|} d\sigma(\By)
  - \frac{\alpha_2}{4\pi} \int_{\p \Omega} \frac{\Bx-\By}{|\Bx-\By|^3} (\Bx-\By) \cdot \mathbf{f}(\By) d\sigma(\By)\\
 & = - \frac{\alpha_1}{4\pi} \int_{\p \Omega} \frac{\mathbf{f}(\By)}{|\Bx-\By|} d\sigma(\By)
  + \frac{\alpha_2}{4\pi} \nabla \int_{\p \Omega} \frac{(\Bx-\By) \cdot \mathbf{f}(\By)}{|\Bx-\By|}  d\sigma(\By) \\
 & \quad\quad  - \frac{\alpha_2}{4\pi} \int_{\p \Omega} \frac{1}{|\Bx-\By|} \nabla_\Bx \left((\Bx-\By) \cdot \mathbf{f}(\By)\right) d\sigma(\By) \\
 & = - \frac{\alpha_1+\alpha_2}{4\pi} \int_{\p \Omega} \frac{\mathbf{f}(\By)}{|\Bx-\By|} d\sigma(\By)
 + \frac{\alpha_2}{4\pi} \nabla \int_{\p \Omega} \frac{(\Bx-\By) \cdot \mathbf{f}(\By)}{|\Bx-\By|}  d\sigma(\By) \\
 & = - \frac{\alpha_1+\alpha_2}{4\pi} \int_{\p \Omega} \frac{\mathbf{f}(\By)}{|\Bx-\By|} d\sigma(\By)
 + \frac{\alpha_2}{4\pi} \nabla  \nabla \cdot \int_{\p \Omega} |\Bx-\By| \mathbf{f}(\By) d\sigma(\By).
 \end{align*}
Since $\Delta |\Bx| = 2|\Bx|^{-1}$, we have
 \beq \label{scalform1}
 \Scal_\Om[\mathbf{f}] (\Bx)
 = - \frac{\alpha_1+\alpha_2}{8\pi} \Delta \int_{\p \Omega} |\Bx-\By| \mathbf{f}(\By) d\sigma(\By)
 + \frac{\alpha_2}{4\pi} \nabla \nabla \cdot \int_{\p \Omega} |\Bx-\By| \mathbf{f}(\By)
 d\sigma(\By) \;.
 \eeq

Let
 \beq
 \Hcal_\Om[\mathbf{f}](\Bx) := \frac{1}{4\pi}\int_{\p \Omega} |\Bx-\By| \mathbf{f}(\By) d\sigma(\By) \;.
 \eeq
Then, in summary, we have
 \beq \label{scalform2}
 \Scal_\Om[\mathbf{f}] (\Bx)
 = - \frac{\alpha_1+\alpha_2}{2} \Delta \Hcal_\Om[\mathbf{f}](\Bx)
 + \alpha_2 \nabla \nabla \cdot \Hcal_\Om[\mathbf{f}](\Bx) \;.
 \eeq
It is worth emphasizing that $\Delta^2 \Hcal_\Om[\mathbf{f}]=0$,
{\it i.e.}, $\Hcal_\Om[\mathbf{f}]$ is biharmonic, in $\Om$ and
$\RR^3 \setminus \overline{\Om}$. Thus \eqnref{scalform2} shows that
the solution to the Lam\'e system in a bounded domain in $\Om$ or
the exterior $\RR^3 \setminus \overline{\Om}$ can be decomposed into
a part harmonic in $\Om$ or $\RR^3 \setminus \overline{\Om}$ and a
gradient part.

Suppose that the solution $\Bu$ to \eqnref{main-eqn} inside $\Om$ is
given by
 \beq\label{inside}
 \Bu (\Bx) = \BB \Bx + \mathbf{v}, \quad \Bx \in \Om
 \eeq
for some constant symmetric matrix $\BB$ and a constant vector
$\mathbf{v}$. Then the solution is given by
 \begin{equation} \label{rep-1}
 \mathbf{u} (\Bx) =
 \begin{cases}
 \BA\Bx + \Scal_\Om [\mathbf{f}] (\Bx)\;, \quad & \Bx \in \RR^3 \setminus
 \overline \Om\;, \\ \BB \Bx + \mathbf{v} \;, \quad & \Bx \in \Om\;,
 \end{cases}
 \end{equation}
where
 \beq
 \mathbf{f}= (\CC_1 - \CC_0) \Ecal(\BB\Bx) \Bn \; .
 \eeq
Here $\Bn=(n_1,n_2,n_3)$ is the unit outward normal vector field to
$\p\Om$. See \cite[Section 4]{KM_ARMA_08}. Note that
 \beq \label{vecscal}
 (\CC_1 - \CC_0) \Ecal(\BB\Bx) \Bn
 = [(\tilde\lambda -\lambda) \mbox{tr} \, (\BB) \mathbf{I} + 2(\tilde\mu-\mu)
 \BB] \Bn .
 \eeq
Let us put
 \beq \label{Bstar}
 \BB^* := (\tilde\lambda -\lambda)
\mbox{tr} \, (\BB) \mathbf{I} + 2(\tilde\mu-\mu) \BB,
 \eeq
so that $\mathbf{f}$ in \eqnref{rep-1} is given by
 \beq\label{fBstar}
 \mathbf{f} = \BB^* \Bn.
 \eeq

According to \eqnref{scalform1}, we have
 \beq
 \Scal_\Om[\BB^*\Bn] (\Bx)
 = - \frac{\alpha_1+\alpha_2}{2} \Delta \Hcal_\Om[\BB^*\Bn](\Bx)
 + \alpha_2 \nabla \nabla \cdot \Hcal_\Om[\BB^*\Bn](\Bx) \;.
 \eeq
One can easily see that
 \beq
 \Hcal_\Om[\BB^*\Bn](\Bx) = \BB^* \Hcal_\Om[\Bn](\Bx)= -\BB^* \nabla p_\Om(\Bx), \quad \Bx \in \Om,
 \eeq
where $p_\Om$ is defined by
 \beq\label{pom}
 p_\Om(\Bx) := \frac{1}{4\pi}\int_{\Omega} |\Bx-\By| d\By \;, \quad \Bx \in
 \RR^3 \; .
 \eeq
Therefore we have
 $$
 \Scal_\Om[\BB^*\Bn] (\Bx) = \frac{\alpha_1+\alpha_2}{2} \BB^* \nabla \Delta p_\Om(\Bx)
 - \alpha_2 \nabla \nabla \cdot \BB^* \nabla p_\Om(\Bx) .
 $$

For a $3 \times 3$ symmetric matrix $\BB$, let $\Delta_\BB:= \nabla
\cdot \BB \nabla$. We note that $\Delta_{\BI} = \Delta$, the usual
Laplacian. We then define
 \beq \label{wom}
 w_\Om^\BB (\Bx) : = \Delta_\BB p_\Om(\Bx), \quad \Bx \in \RR^3.
 \eeq
In particular, we write $w_\Om=w_\Om^{\BI}$. Then, one can easily
see that
 \beq
 w_\Om (\Bx) := \frac{2}{4\pi} \int_{\Om} \frac{1}{|\Bx-\By|}
 d\sigma(\By), \quad \Bx \in \RR^3,
 \label{scal}
 \eeq
which is ($2$ times) the Newtonian potential of $\Om$.

It is appropriate to recall now the proof of the weak Eshelby conjecture
by Kang and Milton. In \cite{KM_ARMA_08}, the matter was reduced
to the statement: `The Newtonian potential is quadratic in $\Om$ if
and only if $\Om$ is an ellipsoid', which was proved by Dive
\cite{div31} and Nikliborc \cite{nik32} in relation to the Newtonian
potential problem (see also \cite{df86}). This statement can be
rephrased as
 \beq\label{dive2}
 \mbox{$w_\Om$ is quadratic in $\Om$ if and only if $\Om$ is an
 ellipsoid,}
 \eeq

If we further put $\alpha:= \frac{\alpha_1+\alpha_2}{2\alpha_2}$,
then we have
 \beq\label{linear2}
 \Scal_\Om[\BB^*\Bn] (\Bx)= \frac{1}{\alpha_2} \left[ \alpha \BB^* \nabla w_\Om(\Bx)
 - \nabla w_\Om^{\BB^*}(\Bx) \right] .
 \eeq
We emphasize that $\alpha > 1$.

%
\section{Proof of Theorem \ref{mainthm}}
%

Suppose that the solution $\Bu$ to \eqnref{main-eqn} is linear in
$\Om$ and given by \eqnref{rep-1}. Then by \eqnref{fBstar} we have
 $$
 \Scal_\Om[\BB^*\Bn] (\Bx) = (\BB-\BA)\Bx + \Bv, \quad \Bx \in \Om.
 $$
It then follows from \eqnref{linear2} that
 \beq\label{linear10}
 \alpha \BB^* \nabla w_\Om(\Bx)
 - \nabla w_\Om^{\BB^*}(\Bx) = \alpha_2 (\BB-\BA)\Bx + \alpha_2 \Bv , \quad \Bx \in \Om.
 \eeq

Note that if eigenvalues of $\BB$ are either all the same or all
distinct, so are those of $\BB^*$. After rotation if necessary, we
may assume that $\BB^*$ is diagonal, say
 \beq
 \BB^* = \mbox{diag} [b_1, b_2, b_3].
 \eeq

(i) Suppose first that all eigenvalues of $\BB^*$ are the same, {\it
i.e.}, $b_1=b_2=b_3=b$. In this case, since $w_\Om^{\BB^*} = b
w_\Om$, it follows from \eqnref{linear10} that
 $$
 b (\alpha-1) \nabla w_\Om = \mbox{linear in } \Om.
 $$
Since $\alpha>1$, $w_\Om$ is quadratic in $\Om$, and hence $\Om$ is
an ellipsoid by \eqnref{dive2}.

 (ii) Suppose now that all eigenvalues of $\BB^*$ are
distinct, {\it i.e.}, $b_i \neq b_j$ if $i \neq j$. In this case,
\eqnref{linear10} yields that
 $$
 \pd{}{x_j} \left(\alpha b_j w_\Om - w_\Om^{\BB^*} \right) = \mbox{linear in }
 \Om, \quad j=1,2,3,
 $$
and hence
 $$
 \alpha b_j w_\Om - w_\Om^{\BB^*} \approx f_j
 (\Bx) \quad \mbox{in } \Om, \quad j=1,2,3,
 $$
for some function $f_j$ which is independent of $x_j$. Here and
afterwards $\approx$ denotes the equality up to a quadratic
function. It then follows that
 \beq\label{aomb}
 \alpha w_\Om \approx \frac{f_1-f_2}{b_1-b_2}, \quad
 w_\Om^{\BB^*} \approx \frac{b_2 f_1- b_1
 f_2}{b_1-b_2},
 \eeq
and
 \beq\label{fonetwo}
 (b_3-b_2) f_1 + (b_1-b_3) f_2 + (b_2-b_1)
 f_3 \approx 0.
 \eeq

Since $f_j$ is independent of $x_j$ for $j=1,2,3$, one can easily
see that \eqnref{fonetwo} holds only when $f_1$, $f_2$ and $f_3$
take the form
 \begin{align*}
 f_1 (\Bx) &\approx \frac{m(x_3)-n(x_2)}{b_3-b_2}, \\
 f_2 (\Bx) &\approx \frac{r(x_1) -m(x_3)}{b_1-b_3}, \\
 f_3 (\Bx) &\approx \frac{n(x_2)-r(x_1)}{b_2-b_1},
 \end{align*}
for some functions $m$, $n$ and $r$. It then follows from
\eqnref{aomb} that
 \begin{align}
 \alpha w_\Om &\approx \frac{m(x_3)}{(b_3-b_2)(b_1-b_3)}
 + \frac{n(x_2)}{(b_2-b_1)(b_3-b_2)}
 + \frac{r(x_1)}{(b_1-b_3)(b_2-b_1)}. \label{awom}
 \end{align}
Since $\Delta w_\Om = 2$ in $\Om$, we have
 $$
 \frac{m''(x_3)}{(b_3-b_2)(b_1-b_3)}
 + \frac{n''(x_2)}{(b_2-b_1)(b_3-b_2)}
 + \frac{r''(x_1)}{(b_1-b_3)(b_2-b_1)}
 = \mbox{constant}.
 $$
Thus $r$, $n$, and $m$ are quadratic functions of $x_1$, $x_2$, and
$x_3$, respectively, and hence $w_\Om$ is quadratic in $\Om$. Thus
$\Om$ is an ellipsoid.

This completes the proof. \qed

%
\section{Proof of Theorem \ref{mainthm2}}
%

In order to prove Theorem \ref{mainthm2}, we use the following lemma
whose proof will be given in the appendix.

\begin{lem}\label{discrim}
Let $\BB_1$ and $\BB_2$ be two symmetric $3\times 3$-matrices. If
$\BB_1+t\BB_2$ has a multiple eigenvalue for all real numbers $t$,
then $\BB_1$ and $\BB_2$ can be diagonalized by the same orthogonal
matrix.
\end{lem}

Let $\BA_1$ and $\BA_2$ be two linearly independent symmetric $3 \times 3$ matrices and
suppose that solutions $\Bu_1$ and $\Bu_2$ to \eqnref{main-eqn} with $\BA=\BA_1$ and $\BA=\BA_2$ are  linear in
$\Om$. Put $\BB_j := \Ecal(\Bu_j)$, $j=1,2$. Since the EMT is positive or negative definite on $M_d^s$ (Theorem 10.6 of \cite{AK_book_07}), \eqnref{EMT2} shows that $\BB_1$ and $\BB_2$ are linearly independent.
According to \eqnref{fBstar} we have
 \begin{align*}
 \Scal_\Om[\BB_1^*\Bn] (\Bx) = (\BB_1-\BA_1)\Bx + \Bv_1, \quad \Bx \in \Om,\\
 \Scal_\Om[\BB_2^*\Bn] (\Bx) = (\BB_2-\BA_2)\Bx + \Bv_2, \quad \Bx \in \Om.
 \end{align*}
It then follows from \eqnref{linear2} that
 \begin{align}\label{linear103}\begin{cases}
 \alpha \BB_1^* \nabla w_\Om(\Bx)
 - \nabla w_\Om^{\BB_1^*}(\Bx) = \alpha_2 (\BB_1-\BA_1)\Bx + \alpha_2 \Bv_1 , \\
 \alpha \BB_2^* \nabla w_\Om(\Bx)
 - \nabla w_\Om^{\BB_2^*}(\Bx) = \alpha_2 (\BB_2-\BA_2)\Bx + \alpha_2 \Bv_2,
 \end{cases}\end{align}
for $ \Bx \in \Om.$

Let us suppose that all of $\BB_1$, $\BB_2$, and $\BB_1+t\BB_2$ $(t
\in \RR)$ have an eigenvalue of multiplicity 2 (otherwise we apply
Theorem \ref{mainthm} to conclude that $\Om$ is an ellipsoid). By
Lemma \ref{discrim}, $\BB_1$ and $\BB_2$ can be diagonalized by a
single orthogonal matrix.  Thus we may assume that $\BB_1$ and
$\BB_2$ are diagonal. Then from \eqnref{Bstar} $\BB_1^*$ and
$\BB_2^*$ are also diagonal and we may let
 $$
 \BB_1^*=\mbox{diag}[b_1,b_1,c_1],\quad \BB_2^*=\mbox{diag}[b_2,b_2,c_2],
 $$
where $b_1\ne c_1$ and $b_2\ne c_2$. Since $\BB_1^*$ and $\BB_2^*$ are linearly independent, we have
 $$
 b_1c_2 \ne c_1 b_2.
 $$

By \eqnref{linear103}, we have
\begin{align}
&\alpha b_1 w_\Om - w^{\BB_1^*}_\Om \approx f(x_3),\label{l1}\\
&\alpha c_1 w_\Om - w^{\BB_1^*}_\Om \approx g(x_1,x_2),\label{l2}\\
&\alpha b_2 w_\Om - w^{\BB_2^*}_\Om \approx h(x_3),\label{l3}\\
&\alpha c_2 w_\Om - w^{\BB_2^*}_\Om \approx l(x_1,x_2),\label{l4}
\end{align}
for some functions $f$, $g$, $h$, and $l$.
Here again $\approx$ denotes the equality up to a quadratic
function. By \eqnref{wom}, we have from \eqnref{l1} and \eqnref{l3}
 \beq
 (\alpha-1)b_2 f(x_3)-(\alpha-1)b_1h(x_3) \approx (\alpha-1)(b_1c_2-b_2c_1) \frac{\partial^2p_\Om}{\partial x_3^2},
 \eeq
and from \eqnref{l2} and \eqnref{l4}
 \beq
 (b_2-\alpha c_2)g(x_1,x_2)-(b_1-\alpha c_1)l(x_1,x_2) \approx
 (1-\alpha)(b_1c_2-b_2c_1) \frac{\partial^2p_\Om}{\partial x_3^2}.
 \eeq
It then follows that
 \begin{equation}\label{p3}
 \frac{\partial^2p_\Om}{\partial x_3^2}  \approx 0.
 \end{equation}
We then obtain from \eqnref{l1}-\eqnref{l4} that
\begin{align*}
(\alpha-1) b_1 w_\Om \approx f(x_3),\quad (\alpha-1) b_2 w_\Om \approx h(x_3),
\end{align*}
and
 $$
 (\alpha c_1-b_1) w_\Om \approx g(x_1,x_2), \quad (\alpha c_2-b_2) w_\Om \approx l(x_1,x_2).
 $$
Thus we conclude that $w_\Om \approx 0$, and hence $\Om$ is an ellipsoid.

This completes the proof.
\qed

%
\section{Proof of Theorem \ref{thm:main-thm}}
%

The space $M_d^s$ is equipped with the inner product $\BA:\BB$,
where $\BA:\BB$ denotes the contraction of two matrices $\BA$ and
$\BB$, {\it i.e.}, $\BA:\BB=\sum_{i,j}
a_{ij}b_{ij}=\textrm{tr}(\BA^{T}\BB)$ where $\textrm{tr}(\BA)$
denotes the trace of $\BA$. For $d=2,3$, let $d_*:=
\frac{d(d+1)}{2}$, which is the dimension of $M_d^s$. Let
$\BB_1=\frac{1}{\sqrt{d}} \mathbf{I}_2$ be a basis for
$\mathbf{\Lambda}_1(M_d^s)$ (of a unit length), and $\{ \BB_2,
\ldots, \BB_{d_*} \}$ be an orthonormal basis for
$\mathbf{\Lambda}_2(M_d^s)$. Then $\{ \BB_1, \ldots, \BB_{d_*} \}$
is an orthonormal basis for $M_d^s$, {\it i.e.},
 $$
 \BB_i : \BB_j = \delta_{ij},
 $$
where $\delta_{ij}$ is Kronecker's delta. Note that for any symmetric
$4$-tensor $\mathbb{T}$, we have
 \beq
 \tr  \mathbb{T} = \sum_{k=1}^{d_*} \mathbb{T} \BB_k : \BB_k \, .
 \eeq

We deal with the case when $\CC^1 - \CC^0$ is positive definite so that $\MM$ is a
symmetric positive-definite linear operator on $M^s_d$. The other
case can be treated in the exactly same way.

Let us first invoke some facts proved in \cite{CK_AMO_08}. Introduce
a $4$-tensor $\widetilde{\MM}$ by
 \begin{equation} \label{wmm}
 \begin{array}{rl}
 \ds \BA:\widetilde{\MM} \BA & = \ds \min_{{\mathbf{v}}\in H^{1}\left(\mathbb{R}^{d}\right)}\int_{\mathbb{R}^{d}}\CC_{\Om}
 \left(\mathcal{E}({\mathbf{v}})+1_{\Om}\GG \BA\right):\left(\mathcal{E}({\mathbf{v}})+1_{\Om}\GG \BA\right)dx\\
 \nm
 &  \ds \quad \quad+\left|\Om\right|\left(\CC^{1}-\CC^{0}\right)\left(\CC^{1}\right)^{-1}\CC^{0}\BA:\BA
 \end{array}
 \end{equation}
for $\BA \in M_d^s$, where
 \begin{equation}
 \GG :=\II-(\CC^{1})^{-1}\CC^{0}.
 \end{equation}
Note that the minimum in \eqnref{wmm} is attained by
$\mathbf{v}=\mathbf{u} - \BA\mathbf{x}$, where $\Bu$ is  the
solution of \eqnref{main-eqn}. It is proved in \cite[Corollary
3.2]{CK_AMO_08} that
 \begin{equation}\label{eq:wmmmm}
  \mathbf{\Lambda}_1 \widetilde{\MM} \mathbf{\Lambda}_1 = \mathbf{\Lambda}_1 \MM \mathbf{\Lambda}_1 \quad \mbox{and} \quad \mathbf{\Lambda}_2 \widetilde{\MM} \mathbf{\Lambda}_2 = \mathbf{\Lambda}_2 \MM \mathbf{\Lambda}_2 \, .
 \end{equation}
In particular, we have
 \begin{equation}
 \tr (\mathbf{\Lambda}_{1}\widetilde{\MM}^{-1}\mathbf{\Lambda}_{1})= \tr (\mathbf{\Lambda}_{1}\MM^{-1}\mathbf{\Lambda}_{1}) \quad \mbox{and} \quad
 \tr (\mathbf{\Lambda}_{2}\widetilde{\MM}^{-1}\mathbf{\Lambda}_{2})=
 \tr (\mathbf{\Lambda}_{2}\MM^{-1}\mathbf{\Lambda}_{2}).
 \end{equation}

Let $\CC$ be an isotropic $4$-tensor, {\it i.e.},
 \[
 \CC=\lambda\mathbf{I}\otimes\mathbf{I} +2\mu\II =d \kappa
 \mathbf{\Lambda_{1}}+2\mu\mathbf{\Lambda_{2}},
 \]
for some $\lambda$ and $\mu$ satisfying \eqnref{eq:trans-2}. Let
$L^{2}\left(\mathbb{R}^{d},M_{d}^{s}\right)$ be the space of square
integrable functions on $\RR^d$ valued in $M_{d}^{s}$ and
$H^1(\mathbb{R}^d, M_d^s)$ the Sobolev space. For $\BP \in
L^{2}\left(\mathbb{R}^{d},M_{d}^{s}\right)$, we define $F_\CC(\BP)$
by
 \begin{equation} \label{fcom}
 F_{\CC}(\BP):= - \mathcal{E} \mathcal{L}_{\CC}^{-1} (\nabla\cdot
 \BP),
 \end{equation}
where $\mathcal{L}_{\CC} = \nabla\cdot \CC \Ecal$. In other words,
if $\Phi$ is a unique solution in $H^1(\mathbb{R}^d, M_d^s)$ to
 \beq\label{LCC}
 \mathcal{L}_{\CC}\left(\Phi \right) + \nabla \cdot \BP = 0,
 \eeq
then $F_{\CC}(\BP)$ is given by
$$
 F_{\CC}(\BP)= \mathcal{E} \left(\Phi \right).
$$

If $\Phi$ is the solution to \eqnref{LCC}, then
$$
\int_{\mathbb{R}^{d}} \CC \left(\mathcal{E} \left(\Phi \right) +
\CC^{-1}\BP\right) :  \mathcal{E} \left(\Psi \right)  =0
$$
for all $\Psi \in H^1(\mathbb{R}^d, M_d^s)$, and hence by taking $\Psi=\Phi$ we have
\begin{equation}\label{eq:FCnegdef}
\int_{\mathbb{R}^{d}} \BP : F_{\CC}(\BP)  = - \int_{\mathbb{R}^{d}}
\CC   F_{\CC}(\BP):F_{\CC}(\BP)
\end{equation}

We prove Theorem~\ref{thm:main-thm} using the following two
propositions whose proofs will be given at the end of this section.

\begin{prop}\label{pro:AC0}
Let
\begin{equation}\label{eq:EAC0}
E_\BA(\CC^0,\BP)= \int_{\Om}
\BP:F_{\CC^0}\left(1_\Om\BP\right)+\int_{\Om}\left(\CC^{0} -
\CC^{1}\right)^{-1}\BP:\BP  + 2\int_{\Om}\BP:\BA \, .
\end{equation}
Then the following holds
\begin{equation}\label{eq:HS-3}
 \BA:\widetilde{\MM} \BA = \sup_{\BP\in L^{2}(\mathbb{R}^{d}:M_{d}^{s}) } E_\BA(\CC^0,\BP).
\end{equation}
Furthermore, this supremum is attained by $\BP=1_\Om
\left(\CC^{1}-\CC^{0}\right) \mathcal{E}(\mathbf{u})$, where
$\mathbf{u}$ is the solution of \eqnref{main-eqn}.
\end{prop}

We then show that structures reaching the lower trace bounds have a
particular structure, as explained by the proposition below.

\begin{prop}\label{pro:trace-opt}
If equality in \eqnref{eq:bd-itrace-1} holds, then we have
 \begin{equation} \label{tropt1}
 E_\BA(\CC^0,1_\Om \BB_1) =
 \sup_{ \BP\in L^{2}(\mathbb{R}^{d}:M_{d}^{s}) }
 E_\BA(\CC^0,\BP)
 \end{equation}
with $\BA= \MM^{-1} \BB_1$. If equality in
\eqnref{eq:bd-itrace-2} holds, then we have
 \begin{equation} \label{tropt2}
 E_\BA(\CC^0,1_\Om \BB_k) =
 \sup_{ \BP\in L^{2}(\mathbb{R}^{d}:M_{d}^{s})}
 E_\BA(\CC^0,\BP)
 \end{equation}
with $\BA= \MM^{-1} \BB_k$, for $k=2, \ldots, d_*$.
\end{prop}

\medskip
\noindent{\sl Proof of Theorem \ref{thm:main-thm}}. Introduce a
bilinear form $\mathcal{F}_{\CC^0}(\BQ,\BR)$ by
$$
\mathcal{F}_{\CC^0}(\BQ,\BR) =  \int_{\Om}
\BQ:F_{\CC^0}\left(1_\Om\BR\right) +\int_{\Om}\left(\CC^{0} -
\CC^{1}\right)^{-1}\BQ:\BR \, .
$$
It follows from \eqnref{eq:FCnegdef} that
 \begin{align*}
 \mathcal{F}_{\CC^0}(1_\Om \BQ,1_\Om \BQ) & =
 -  \int_{\Om} \CC^0 F_{\CC^0}\left(1_\Om\BQ\right) :F_{\CC^0}\left(1_\Om\BQ\right)
 - \int_{\Om}\left(\CC^{1} - \CC^{0}\right)^{-1}\left(1_\Om\BQ\right):\left(1_\Om\BQ\right) \\
 &\leq - \int_{\Om}\left(\CC^{1} - \CC^{0}\right)^{-1}\left(1_\Om\BQ\right):\left(1_\Om\BQ\right) \\
 &\leq - K \left\| 1_\Om\BQ\right\|_{L^{2}(\mathbb{R}^{d}:M_{d}^{s})}
 \end{align*}
for some positive constant $K$. The last holds due to the
positive-definiteness of $\CC^{1} - \CC^{0}$. As a consequence,
$\mathcal{F}_{\CC^0}$ is negative definite when restricted to
$H=\{\BP \in L^{2}(\Om:M_{d}^{s})^2, \mbox{ supported in } \Om \}$.
Therefore $E_\BA(\CC^0,\BQ)=\mathcal{F}_{\CC^0}(\BQ,\BQ) +
2\int_{\Om}\BQ:\BA $ is a strictly concave functional on $H$, and
therefore admits at most one maximizer in $H$.

We observe that since $\CC^{1}-\CC^{0}$ is isotropic, if $\BB$ is
diagonal and all the eigenvalues are the same, so is
$(\CC^{1}-\CC^{0})^{-1} \BB$. If $\BB$ is trace-free and all the
eigenvalues are distinct, so is $(\CC^{1}-\CC^{0})^{-1} \BB$.

Suppose that equality holds in \eqnref{eq:bd-itrace-1}. It then
follows from Proposition~\ref{pro:AC0}, \eqnref{pro:trace-opt}, and
uniqueness of the maximizer in $\Om$ that
$$
\left(\CC^{1}-\CC^{0}\right)\mathcal{E}(\mathbf{u}_1) = \BB_1
\textrm{ in } \Om,
$$
where $\Bu_1$ is the solution to \eqnref{main-eqn} with
$\BA=\MM^{-1} \BB_1$. Recall that $\BB_1=\frac{1}{\sqrt{d}} \BI$.
Therefore, $\Ecal(\Bu_1)$ is constant in $\Om$ and all the
eigenvalues of $\Ecal(\Bu_1)$ are the same. Thus $\Om$ is an ellipse
or an ellipsoid due to Theorem \ref{mainthm}.

Suppose now that equality holds in \eqnref{eq:bd-itrace-2}. Then
for similar reasons we can deduce that for each $k=2, \ldots, d_*$,
$$
\left(\CC^{1}-\CC^{0}\right)\mathcal{E}(\mathbf{u}_k) = \BB_k
\textrm{ in } \Om,
$$
where $\Bu_k$ is the solution to \eqnref{main-eqn} with
$\BA=\MM^{-1} \BB_k$. Thus $\Om$ is an ellipse or an ellipsoid
due to Theorem \ref{mainthm2}.

This completes the proof. \qed

\medskip

\noindent{\sl Proof of Proposition \ref{pro:AC0}}. Following the
notation of \cite{CK_AMO_08}, we  define $W_\BA(\CC,\BP)$, for $\BA
\in M_d^s$, by
 $$
 W_\BA(\CC,\BP)=
 \int_{\mathbb{R}^{d}}\BP:F_{\CC}\BP+\int_{\mathbb{R}^{d}}\left(\CC- \CC_{\Om}\right)^{-1}\BP:\BP
 + 2\int_{\Om}\!\BP:\left(\CC^{1}-\CC\right)^{-1}\left(\CC^{1}-\CC^{0}\right)\BA .
 $$
It is proved in \cite[Proposition 4.1]{CK_AMO_08}, following the variational strategy given in \cite{KM_IMA_86} for the derivation of Hashin-Shtrikman type bounds, that for any isotropic elasticity tensor $\CC<\CC^{0} (< \CC^1)$ we have
 \begin{equation}\label{eq:HS-1}
 \BA:\widetilde{\MM} \BA= \BA:\left(\CC^{1}-\CC^{0}\right)\left(\CC-\CC^{1}\right)^{-1}\left(\CC-\CC^{0}\right)\BA
 + \sup_{\BP\in L^{2}\left(\mathbb{R}^{d}:M_{d}^{s}\right)}W_\BA(\CC,\BP).
 \end{equation}
Note that the supremum is attained by
\begin{equation}\label{eq:sigmaopt}
\BP = 1_\Om \left(\CC^{1}-\CC^{0}\right) \BA +
\left(\CC_\Om-\CC\right)\left(\mathcal{E}(\mathbf{u})-\BA\right),
\end{equation}
where $\Bu$ is the solution to \eqnref{main-eqn} with $\BA$ in above
identity. Since $(\CC^{1}-\CC^{0}) (\CC-\CC^{1})^{-1} (\CC-\CC^{0})$
is positive definite, by sending $\CC$ to $\CC^0$, and restricting
the supremum to fields $\BP$ such that $1_\Om \BP =\BP$ we obtain
\begin{equation}
 \BA:\widetilde{\MM} \BA \ge \sup_{ \BP\in L^{2}(\mathbb{R}^{d},M_{d}^{s})}
 E_\BA(\CC^0,\BP) \, .
 \end{equation}

For any $\BP \in  L^{2}\left(\mathbb{R}^{d},M_{d}^{s}\right)$, and
any positive definite isotropic elasticity tensor $\CC < \CC^0$, we
define $E_\BA(\CC,\BP)$, for $\BA \in M_d^s$, by
 $$
 E_\BA(\CC,\BP)= \int_{\Om}\! \BP:\!F_{\CC}\left(1_\Om\BP\right)+\int_{\Om}\!\left(\CC - \CC^{1}\right)^{-1}\!\BP:\BP 
 + 2\int_{\Om}\!\BP:\!\left(\CC^{1}-\CC\right)^{-1}\!\left(\CC^{1}-\CC^{0}\right)\!\BA .
 $$
Note that this definition is consistent of that of
$E_\BA(\CC^0,\BP)$ given in \eqnref{eq:EAC0} by passing to the limit
in $\CC$.

Introducing the decomposition $\BP = \BP_\Om +\BP_U$, with $\BP_\Om
1_\Om =\BP_\Om$, and $\BP_U 1_\Om \equiv 0$, we have
\begin{eqnarray*}
 W_\BA(\CC,\BP)&=&E_\BA(\CC,\BP_\Om) + W_\BA(\CC,\BP_U) \\
&+& \int_{\mathbb{R}^{d}}\BP_\Om:F_{\CC}\BP_U +
\int_{\mathbb{R}^{d}}\BP_U:F_{\CC}\BP_\Om.
\end{eqnarray*}

Let $\CC=\CC_0 - \eps \II$, where $\eps>0$. Then we have
$$
W_\BA(\CC,\BP)= E_\BA(\CC,\BP_\Om) -
\epsilon^{-1}\|\BP_U\|^2_{L^2\left(\mathbb{R}^{d}\right)} +
R(\BP_U,\BP_\Om),
$$
where
$$
R\left(\BP_U,\BP_\Om\right)   :=
\int_{\mathbb{R}^{d}}\BP_\Om:F_{\CC}\BP_U   +
\int_{\mathbb{R}^{d}}\BP_U:F_{\CC}\BP_\Om
   + \int_{\mathbb{R}^{d}}\BP_U:F_{\CC}\BP_U .
$$
By integration by parts, and by the Cauchy-Schwartz inequality, we
readily obtain that for $\eps$ small enough
$R\left(\BP_U,\BP_\Om\right)$  satisfies
$$
\left|R\left(\BP_U,\BP_\Om\right) \right| \leq  K
\|\BP_U\|_{L^2\left(\mathbb{R}^{d}\right)}
       \left( \|\BP_U\|_{L^2\left(\mathbb{R}^{d}\right)}
            + \|\BP_\Om\|_{L^2\left(\mathbb{R}^{d}\right)}
       \right),
$$
where the constant $K$ is independent of $\BP_U$, $\BP_\Om$, and
$\eps$. As a consequence, for $\eps$ small enough,
$$
- \epsilon^{-1}\|\BP_U\|^2_{L^2\left(\mathbb{R}^{d}\right)} +
R(\BP_U,\BP_\Om) \leq 3K \eps
\|\BP_\Om\|_{L^2\left(\mathbb{R}^{d}\right)}^2.
$$
Note that from \eqnref{eq:FCnegdef} $\int_{\Om}
\BP_\Om:F_{\CC}\left(\BP_\Om\right)$ is negative definite, therefore
$$
E_{\BA}(\CC,\BP)\leq \tilde{K}
\|\BP_\Om\|_{L^2\left(\mathbb{R}^{d}\right)} (-
\|\BP_\Om\|_{L^2\left(\mathbb{R}^{d}\right)} + 1),
$$
where $\tilde{K}$ is another constant independent of $\BP_\Om$ and
$\eps$. Thus $\|\BP_\Om\|_{L^2\left(\mathbb{R}^{d}\right)}$ must
stay bounded, uniformly with respect to $\epsilon$, close to the
supremum.  Taking the limit as $\eps$ tends to zero we obtain
\eqnref{eq:HS-3}. Replacing $\CC$ by $\CC^0$ in \eqnref{eq:sigmaopt}
concludes the proof. \qed

\medskip

\noindent{\sl Proof of Proposition~\ref{pro:trace-opt}}. Given
$k\in\{1,\ldots,d_*\}$, choose $\BA= \MM^{-1}
\mathbf{\Lambda}_l(\BB_k)$ and use a test function $\BP = 1_\Om
\mathbf{\Lambda}_l(\BB_k)$ in \eqnref{eq:HS-3}. This gives
 \begin{align}
 & \MM^{-1} \mathbf{\Lambda}_l(\BB_k): \mathbf{\Lambda}_l(\BB_k) \ge W_\BA(\CC^0, 1_\Om
 \mathbf{\Lambda}_l(\BB_k)) \label{eq:mm-one} \\
 &\quad = \int_{\Om} \mathbf{\Lambda}_l(\BB_k) : F_{\CC^0} \left( 1_\Om \mathbf{\Lambda}_l(\BB_k) \right)
 +\int_{\Om}\left(\CC^0- \CC^1 \right)^{-1} \mathbf{\Lambda}_l(\BB_k) :
 \mathbf{\Lambda}_l(\BB_k) \nonumber \\
 &\quad \quad + 2\int_{\Om}\mathbf{\Lambda}_l(\BB_k) : \MM^{-1} \mathbf{\Lambda}_l(\BB_k)  \, . \nonumber
 \end{align}
Summing these inequalities over $k$, we obtain
 \begin{align}
 \tr  (\mathbf{\Lambda}_l \MM^{-1} \mathbf{\Lambda}_l) & \ge \sum_{k=1}^{d_*}
 \int_{\Om} \mathbf{\Lambda}_l(\BB_k) : F_{\CC^0} \left( 1_\Om \mathbf{\Lambda}_l(\BB_k)
 \right) \label{eq:trlam} \\
 & \quad + \sum_{k=1}^{d_*} \int_{\Om}\left(\CC^0- \CC^1 \right)^{-1} \mathbf{\Lambda}_l(\BB_k) :
 \mathbf{\Lambda}_l(\BB_k) + 2\tr  (\mathbf{\Lambda}_l \MM^{-1} \mathbf{\Lambda}_l)  \, .
 \nonumber
 \end{align}
It is proved in \cite[(4.27) \& (4.28)]{CK_AMO_08} that
 $$
 \sum_{k=1}^{d_*}
 \int_{\Om} \mathbf{\Lambda}_1(\BB_k) : F_{\CC^0} \left( 1_\Om \mathbf{\Lambda}_1(\BB_k)
 \right) = - \frac{1}{d\left(\lambda +2\mu \right)} \, ,
 $$
and
 $$
 \sum_{k=1}^{d_*}
 \int_{\Om} \mathbf{\Lambda}_2(\BB_k) : F_{\CC^0} \left( 1_\Om \mathbf{\Lambda}_2(\BB_k)
 \right)
 = - \left(\frac{d-1}{d\left(\lambda +2\mu \right)}+\frac{d-1}{2\mu} \right) \, .
 $$
Since
 $$
 \left(\CC^0- \CC^1 \right)^{-1} = \frac{1}{d(\kappa-\tilde\kappa)}
 \mathbf{\Lambda}_1 + \frac{1}{2(\mu-\tilde\mu)} \mathbf{\Lambda}_2 \, ,
 $$
one can immediately see that
 $$
 \sum_{k=1}^{d_*} \int_{\Om}\left(\CC^0- \CC^1 \right)^{-1} \mathbf{\Lambda}_1(\BB_k) :
 \mathbf{\Lambda}_1(\BB_k) = \frac{1}{d(\kappa-\tilde\kappa)}
 $$
and
 $$
 \sum_{k=1}^{d_*} \int_{\Om}\left(\CC^0- \CC^1 \right)^{-1} \mathbf{\Lambda}_2(\BB_k) :
 \mathbf{\Lambda}_2(\BB_k) = \frac{d_*-1}{2(\mu-\tilde\mu)} \, .
 $$
Therefore, we get
 $$
 \sum_{k=1}^{d_*} \left[ \int_{\Om} \mathbf{\Lambda}_l(\BB_k) : F_{\CC^0} \left( 1_\Om \mathbf{\Lambda}_l(\BB_k)
 \right)+ \int_{\Om}\left(\CC^0- \CC^1 \right)^{-1} \mathbf{\Lambda}_l(\BB_k) :
 \mathbf{\Lambda}_l(\BB_k) \right]  = - K_l
 $$
for $l=1,2$, where $K_l$ is given in \eqnref{eq:kone} and
\eqnref{eq:ktwo}. It then follows from \eqnref{eq:trlam} that

\begin{equation}\label{eq:crucial}
 \tr  (\mathbf{\Lambda}_l \MM^{-1} \mathbf{\Lambda}_l) \ge -K_l + 2\tr  (\mathbf{\Lambda}_l \MM^{-1} \mathbf{\Lambda}_l).
\end{equation}

Suppose that equality in \eqnref{eq:bd-itrace-1} holds. Then, in
view of \eqnref{eq:crucial}, the inequality in \eqnref{eq:trlam}
becomes an equality, and so does the one in \eqnref{eq:mm-one}.
Since $\mathbf{\Lambda}_1(\BB_1)=\BB_1$ and
$\mathbf{\Lambda}_1(\BB_k)=0$ for $k=2, \ldots, d_*$, we have

 \begin{equation}
 E_\BA(\CC^0,1_\Om \BB_1) =
 \sup_{ \BP\in L^{2}(\mathbb{R}^{d}:M_{d}^{s}) }
 E_\BA(\CC^0,\BP) \, , \quad \BA= \MM^{-1} \BB_1.
 \end{equation}
Likewise, if equality in \eqnref{eq:bd-itrace-2} holds, then
 \begin{equation}
 E_\BA(\CC^0,1_\Om \BB_k) =
 \sup_{ \BP\in L^{2}(\mathbb{R}^{d}:M_{d}^{s})}
 E_\BA(\CC^0,\BP) \, , \quad \BA= \MM^{-1} \BB_k,
 \end{equation}
for $k=2, \ldots, d_*$, and the proof is complete. \qed

\appendix
\section{Proof of Lemma \ref{discrim}}

\pf By considering $\frac{1}{t} \BB_1 + \BB_2$ and taking the limit
$t \to \infty$, one can see that $\BB_2$ also has a multiple
eigenvalue. If $\BB_2$ has an eigenvalue of multiplicity 3, then
$\BB_2$ is a constant multiple of the identity matrix and hence the
conclusion of the lemma holds trivially.

Let us assume that $\BB_2$ has an eigenvalue of multiplicity 2. Note
that, for any real number $s$ and orthogonal matrix $\BU$,
 $$
 \BU\BB_1\BU^{-1}+t\BU\BB_2\BU^{-1}- s\BI = \BU (\BB_1 +t\BB_2- s\BI) \BU^{-1}
 $$
has a multiple eigenvalue regardless of  $t$. Therefore we may
assume that $\BB_2$ takes the form
 $$
 \BB_2=\begin{pmatrix} 0 & 0 & 0 \\ 0 & 0 & 0 \\ 0 & 0 & 1 \end{pmatrix}
 $$
while $\BB_1$ is arbitrary, say
 $$
 \BB_1=\begin{pmatrix} a & d & e \\ d & b & f \\ e & f & c \end{pmatrix}.
 $$

Let $\Gamma(t)$ be the discriminant of the characteristic polynomial
of $\BB_1+t\BB_2$. Since $\BB_1+t\BB_2$ has a multiple eigenvalue
for all $t$, $\Gamma(t) \equiv 0$. Then a straightforward
calculation shows that the coefficient of $t^4$ term of $\Gamma(t)$
is given by $a^2+b^2-2ab+4d^2,$ and hence we have
 $$ a-b=d=0.$$
It then follows that the coefficient of $t^2$ term of $\Gamma(t)$ is
given by $(e^2+f^2)^2$, and hence
 $$ e=f=0. $$
Thus $\BB_1$ takes the form
 $$
 \BB_1=\begin{pmatrix} a & 0 & 0 \\ 0 & a & 0 \\ 0 & 0 & c \end{pmatrix}.
 $$

This completes the proof. \qed

\end{document}